# Surprise probabilities in Markov chains

James Norris [*]
University of Cambridge

Yuval Peres
Microsoft Research

Alex Zhai
Stanford University

August 6, 2014

**Abstract**

In a Markov chain started at a state $x$, the *hitting time* $\tau(y)$ is the first time that the chain reaches another state $y$. We study the probability $\mathbf{P}_x(\tau(y) = t)$ that the first visit to $y$ occurs precisely at a given time $t$. Informally speaking, the event that a new state is visited at a large time $t$ may be considered a "surprise". We prove the following three bounds:

- In any Markov chain with $n$ states, $\mathbf{P}_x(\tau(y) = t) \leq \frac{n}{t}$.
- In a reversible chain with $n$ states, $\mathbf{P}_x(\tau(y) = t) \leq \frac{\sqrt{2n}}{t}$ for $t \geq 4n + 4$.
- For random walk on a simple graph with $n \geq 2$ vertices, $\mathbf{P}_x(\tau(y) = t) \leq \frac{4e \log n}{t}$.

We construct examples showing that these bounds are close to optimal. The main feature of our bounds is that they require very little knowledge of the structure of the Markov chain.

To prove the bound for random walk on graphs, we establish the following estimate conjectured by Aldous, Ding and Oveis-Gharan (private communication): For random walk on an $n$-vertex graph, for every initial vertex $x$,

$$\sum_y \left( \sup_{t \geq 0} p^t(x, y) \right) = O(\log n).$$

[*]partially supported by EPSRC grant EP/I03372X/1

# 1 Introduction

Suppose that a Markov chain with $n$ states is run for a long time $t \gg n$. It would be surprising if the state visited on the $t$-th step was not visited at any earlier time — if a state is likely, then we expect it to have been visited earlier, and if it is rare, we do not expect that its first visit will occur at precisely time $t$. How surprised should we be?

Let $X = \{X_t\}$ be a Markov chain with finite state space $S$ of size $n$, and let $\mathbf{P}_x$ denote the probability for the chain started at $x$. Let $\mathcal{S}_t$ denote the event that the state visited at time $t$ was not seen at any previous time. In this paper, we quantify the intuition that $\mathcal{S}_t$ is unlikely by proving upper bounds on $\mathbf{P}(\mathcal{S}_t)$. This question was posed by A. Kontorovich (private communication), who asked for bounds that do not require detailed knowledge of the transition probabilities of the chain.

For $y \in S$, the *hitting time* $\tau(y)$ is defined as the first time that the chain reaches $y$. We can express $\mathcal{S}_t$ in terms of the hitting time as follows.

$$\mathcal{S}_t = \bigcup_{y \in S} \{\tau(y) = t\}.$$

If we can prove a statement of the form $\mathbf{P}_x(\tau(y) = t) \leq M$ for all $y \in S$, then it follows that $\mathbf{P}_x(\mathcal{S}_t) \leq n \cdot M$. Although naive, it turns out that this approach gives bounds on $\mathbf{P}_x(\mathcal{S}_t)$ which are close to optimal, as we will explain in more detail later.

We start with a simple proposition which bounds $\mathbf{P}_x(\tau(y) = t)$ in the most general setting.

**Proposition 1.1.** *Let $X$ be a Markov chain with finite state space $S$ of size $n$, and let $x$ and $y$ be any two states. Then, for all $t > n$,*

$$\mathbf{P}_x(\tau(y) = t) \leq \frac{n}{t}, \quad \text{whence} \quad \mathbf{P}_x(\mathcal{S}_t) \leq \frac{n^2}{t}.$$

Our first main theorem improves this bound for reversible chains.

**Theorem 1.2.** *Let $X$ be a reversible Markov chain with finite state space $S$ of size $n$, and let $\pi$ be the stationary distribution of $X$. Consider any two states $x$ and $y$. For all $t > 0$,*

$$\mathbf{P}_x(\tau(y) = t) \leq \frac{2e \cdot \max\left(1, \log \frac{1}{\pi(x)}\right)}{t}.$$

*In particular, if $X$ is the random walk on a simple graph with $n \geq 2$ vertices, then*

$$\mathbf{P}_x(\tau(y) = t) \leq \frac{4e \cdot \log n}{t} \quad \text{and} \quad \mathbf{P}_x(\mathcal{S}_t) \leq \frac{4e \cdot n \log n}{t}.$$

One limitation of Theorem 1.2 is that in general, the stationary probability $\pi(x)$ can be arbitrarily small. Our second main theorem gives an alternate bound for reversible chains that, in the style of Proposition 1.1, depends only on $n$ and $t$.

**Theorem 1.3.** *Let $X$ be a finite reversible Markov chain with $n$ states, and let $x$ and $y$ be any two states. Then, for all $t \geq 4n + 4$,*

$$\mathbf{P}_x(\tau(y) = t) \leq \frac{\sqrt{2n}}{t}, \quad \text{whence} \quad \mathbf{P}_x(\mathcal{S}_t) \leq \frac{n\sqrt{2n}}{t}.$$

**Remark 1.4.** For reversible chains with non-negative eigenvalues, we prove the stronger bound $\mathbf{P}_x(\tau(y) = t) \leq \frac{1}{2}\sqrt{\frac{n}{t(t-n)}}$. The degrading of the estimate as $t$ approaches $n$ cannot be avoided; if $t < n$, then it is possible to have $\mathbf{P}_x(\tau(y) = t)$ arbitrarily close to 1 by considering a birth-and-death chain on states $1, 2, \ldots, n$ where state $i$ transitions to state $i+1$ with probability $1 - \epsilon$ for $\epsilon \to 0$.

To give a sense of how accurate these estimates are, we give several constructions where $\mathbf{P}_x(\tau(y) = t)$ can be relatively large. In particular, we will show that Proposition 1.1 and Theorem 1.3 are tight up to a constant factor. We also construct a family of simple graphs which achieve

$$\mathbf{P}_x(\tau(y) = t) \geq \frac{c\sqrt{\log n}}{t} \tag{1}$$

for some constant $c > 0$. This does not quite match the upper bound of $(4e \log n)/t$ appearing in Theorem 1.2, but it demonstrates that the dependence on $n$ cannot be avoided. All of these constructions can be modified slightly to give lower bounds on $\mathbf{P}_x(\mathcal{S}_t)$ which are on the order of $n$ times the corresponding lower bounds on $\mathbf{P}_x(\tau(y) = t)$.

In the proof of Theorem 1.2, we show a certain "maximal probability" bound that may be of independent interest. For any two states $x, y \in S$, define

$$p^*(x, y) = \sup_{t \geq 0} p^t(x, y),$$

where $p^t(x, y)$ is the transition probability from $x$ to $y$ in $t$ steps. It was asked by Aldous, and independently by Ding and Oveis Gharan, whether for random walks on simple graphs with $n$ vertices, $\sum_{y \in S} p^*(x, y) = O(\log n)$ for any starting vertex $x$ (private communication). Using a theorem of Starr [9], we prove the following proposition, which verifies this.

**Proposition 1.5.** *Let $X$ be a reversible Markov chain with finite state space $S$ and stationary distribution $\pi$. Then, for any $x \in S$,*

$$\sum_{y \in S} p^*(x, y) \leq 2e \cdot \max\left(1, \log \frac{1}{\pi(x)}\right).$$

**Remark 1.6.** For simple random walk on an $n$-vertex graph, the right hand side is at most $4e \cdot \log n$. When the graph is a cycle, this bound is tight up to a constant factor (see the end of Section 2.1).

Finally, we mention some situations where stronger bounds for $\mathbf{P}_x(\tau(y) = t)$ on the order of $\frac{1}{t}$ are possible.

**Proposition 1.7.** *Let $X$ be a Markov chain with state space $S$ and stationary distribution $\pi$. Then, for any $y \in S$ and $t > 0$,*

$$\mathbf{P}_\pi\left(\tau(y) = t\right) \leq \frac{1}{t},$$

*where $\mathbf{P}_\pi$ denotes the probability for the walk started from a random state sampled from $\pi$. In particular, for any state $x$,*

$$\mathbf{P}_x(\tau(y) = t) \leq \frac{1}{t\pi(x)}.$$

A consequence of Proposition 1.7 is that $\mathbf{P}_x(\tau(y) = t) = O\left(\frac{1}{t}\right)$ when $t$ is significantly larger than the mixing time. We define

$$d_x(t) = \left\|p^t(x, \cdot) - \pi\right\|_{\mathrm{TV}},$$

where $\|\cdot\|_{\mathrm{TV}}$ denotes total variation distance, and recall that the mixing time $t_{\mathrm{mix}}(\epsilon)$ is defined as the earliest time $t$ for which $d_x(t) \leq \epsilon$ for every state $x$.

**Proposition 1.8.** *For a Markov chain $X$ with state space $S$ of size $n$, suppose that we have a bound of the form $\mathbf{P}_x(\tau(y) = t) \leq \psi(t)$ for all $x, y \in S$ (e.g., the bounds of Proposition 1.1 or Theorems 1.2 and 1.3). For any $t > s > 0$,*

$$\mathbf{P}_x(\tau(y) = t) \leq d_x(s)\psi(t - s) + \frac{1}{t - s}.$$

*In particular, if $t > 2t_{mix}(1/4)\lceil \log_2 n \rceil$, we have*

$$\mathbf{P}_x(\tau(y) = t) \leq \frac{4}{t}.$$

The rest of the paper is organized as follows. In Section 2, we prove Proposition 1.1 and Theorem 1.2, which are proved by the same method. We prove Proposition 1.5 as an intermediate step, and we also describe constructions giving lower bounds on $\mathbf{P}_x(\tau(y) = t)$.

Section 3 is self-contained; we record bounds on sums of independent geometric random variables. These are used in Section 4 to prove Theorem 1.3. Here, we also construct an example to show that the bound is of the right order.

Section 5 is devoted to describing the modified constructions which give lower bounds for $\mathbf{P}_x(\mathcal{S}_t)$. We prove Propositions 1.7 and 1.8 in Section 6 and discuss open problems in Section 7. The last two sections contain acknowledgements and background.

## 2 Proofs of Proposition 1.1 and Theorem 1.2

We start with a bound on $\mathbf{P}_x(\tau(y) = t)$ in terms of the maximal probabilities $p^*(x, y)$.

**Lemma 2.1.** *Let $X$ be a Markov chain with finite state space $S$, and let $x$ and $y$ be any two states. Then, for all $t > 0$,*

$$\mathbf{P}_x(\tau(y) = t) \leq \frac{1}{t} \sum_{z \in S} p^*(x, z).$$

*Proof.* Observe that for each time $s < t$,

$$\mathbf{P}_x(\tau(y) = t) \leq \sum_{z \in S} p^s(x, z) \mathbf{P}_z(\tau(y) = t - s).$$

Summing this inequality over all $s = 0, 1, 2, \ldots, t - 1$, we obtain

$$t \cdot \mathbf{P}_x(\tau(y) = t) \leq \sum_{s=0}^{t-1} \sum_{z \in S} p^s(x, z) \mathbf{P}_z(\tau(y) = t - s)$$

$$= \sum_{z \in S} \sum_{s=0}^{t-1} p^s(x, z) \mathbf{P}_z(\tau(y) = t - s) \leq \sum_{z \in S} p^*(x, z) \sum_{s=0}^{t-1} \mathbf{P}_z(\tau(y) = t - s)$$

$$= \sum_{z \in S} p^*(x, z) \mathbf{P}_z(1 \leq \tau(y) \leq t) \leq \sum_{z \in S} p^*(x, z).$$

Dividing both sides by $t$ completes the proof. □

*Proof of Proposition 1.1.* By Lemma 2.1, we have

$$\mathbf{P}_x(\tau(y) = t) \leq \frac{1}{t} \sum_{z \in S} p^*(x, z) \leq \frac{1}{t} \sum_{z \in S} 1 = \frac{n}{t}.$$

□

Under the assumption that Proposition 1.5 holds, we can now also prove Theorem 1.2.

*Proof of Theorem 1.2.* By Lemma 2.1 and Proposition 1.5,

$$\mathbf{P}_x(\tau(y) = t) \leq \frac{1}{t}\sum_{z \in S} p^*(x,z) \leq \frac{2e \cdot \max\left(1, \log \frac{1}{\pi(x)}\right)}{t}.$$

For the random walk on a simple graph having $n \geq 2$ vertices, the stationary probability of a vertex is proportional to its degree, and in particular it is at least $\frac{1}{n^2}$. Thus, in this case

$$\mathbf{P}_x(\tau(y) = t) \leq \frac{4e \cdot \log n}{t}.$$

$\square$

## 2.1 Proof of Proposition 1.5

It remains to prove Proposition 1.5. For a Markov chain $X = \{X_t\}$, let $P$ be the transition operator of the chain. That is, for a function $f : S \to \mathbb{R}$, we have

$$(Pf)(x) = \sum_{y \in S} p(x,y)f(y).$$

We deduce Proposition 1.5 from the following theorem of Starr [9].

**Theorem 2.2** (special case of [9], Theorem 1). *Let $X = \{X_t\}$ be a reversible Markov chain with state space $S$ and stationary measure $\pi$. Then, for any $1 < p < \infty$ and $f \in L^p(S, \pi)$,*

$$\mathbf{E}_\pi \left|\sup_{n \geq 0} P^{2n} f\right|^p \leq \left(\frac{p}{p-1}\right)^p \mathbf{E}_\pi |f|^p,$$

*where $\mathbf{E}_\pi$ denotes expectation with respect to $\pi$.*

We include a short proof of Theorem 2.2 in Section 9.1.

*Proof of Proposition 1.5.* Define the two quantities

$$p_e^*(x,y) = \sup_{t \geq 0} p^{2t}(x,y) \qquad \text{and} \qquad p_o^*(x,y) = \sup_{t \geq 0} p^{2t+1}(x,y).$$

for even and odd times. Take $f(y) = \delta_{x,y}$ in Theorem 2.2. Then, we find that for each $t$,

$$(P^t f)(y) = \sum_{z \in S} p^t(y,z) f(z) = p^t(y,x) = \frac{p^t(x,y)\pi(x)}{\pi(y)},$$

and so

$$\sup_{t \geq 0}(P^{2t}f)(y) = \frac{p_e^*(x,y)\pi(x)}{\pi(y)}.$$

We now apply Theorem 2.2 with an exponent $p$ to be specified later. We obtain

$$\sum_{y \in S} \pi(y)\left(\frac{p_e^*(x,y)\pi(x)}{\pi(y)}\right)^p \leq \left(\frac{p}{p-1}\right)^p \cdot \pi(x), \qquad (2)$$

Let $q = \frac{p}{p-1}$ be the conjugate exponent of $p$, so that $\frac{1}{p} + \frac{1}{q} = 1$. By Hölder's inequality,

$$\sum_{y \in S} p_e^*(x,y) = \sum_{y \in S} \pi(y)^{\frac{1}{q}} \frac{p_e^*(x,y)}{\pi(y)^{\frac{1}{q}}} \leq \left( \sum_{y \in S} \pi(y) \right)^{\frac{1}{q}} \left( \sum_{y \in S} \left( \frac{p_e^*(x,y)}{\pi(y)^{\frac{1}{q}}} \right)^p \right)^{\frac{1}{p}}$$

$$= \frac{1}{\pi(x)} \left( \sum_{y \in S} \pi(y) \left( \frac{p_e^*(x,y) \pi(x)}{\pi(y)} \right)^p \right)^{\frac{1}{p}} \leq \frac{p}{p-1} \cdot \pi(x)^{\frac{1}{p}-1} = q \pi(x)^{-\frac{1}{q}},$$

where the second inequality comes from (2). This is valid for all $1 < q < \infty$, and by continuity, also for $q = 1$. Setting $q = \max\left(1, \log \frac{1}{\pi(x)}\right)$, this yields

$$\sum_{y \in S} p_e^*(x,y) \leq e \cdot \max\left(1, \log \frac{1}{\pi(x)}\right).$$

We also have

$$\sum_{y \in S} p_o^*(x,y) \leq \sum_{y \in S} \sum_{z \in S} p_e^*(x,z) p(z,y) = \sum_{z \in S} p_e^*(x,z) \sum_{y \in S} p(z,y)$$

$$= \sum_{z \in S} p_e^*(x,z) \leq e \cdot \max\left(1, \log \frac{1}{\pi(x)}\right),$$

and so

$$\sum_{y \in S} p^*(x,y) \leq \sum_{y \in S} p_e^*(x,y) + \sum_{y \in S} p_o^*(x,y) \leq 2e \cdot \max\left(1, \log \frac{1}{\pi(x)}\right),$$

which proves the lemma. $\square$

We mention an example that shows Proposition 1.5 is tight up to a multiplicative constant. Consider a cycle of even size $n = 2m$, with vertices labeled by elements of $\mathbb{Z}/n$ in the natural way. Let $X$ be the simple random walk started at 0. Then, for each time $t \geq 0$ and each position $-m < k \leq m-1$, the probability $p^t(0,k)$ that $X$ is at $k$ is at least the probability $\tilde{p}^t(0,k)$ that a simple random walk on $\mathbb{Z}$ started at 0 is at $k$ after time $t$. By the local central limit theorem (see e.g. [5], Theorem 1.2.1), we have for some universal constant $C > 0$ that

$$p^t(0,k) \geq \tilde{p}^t(0,k) \geq \frac{2}{\sqrt{2\pi t}} \exp\left(-\frac{k^2}{2t}\right) - \frac{C}{\sqrt{t} \cdot k^2}.$$

Plugging in $t = k^2$ (which approximately maximizes the right hand side), we obtain

$$p^*(0,k) \geq \frac{2}{\sqrt{2\pi e}} \cdot \frac{1}{|k|} - \frac{C}{k^3}.$$

Thus,

$$\sum_{k \in \mathbb{Z}/n} p^*(0,k) \geq \sum_{k=1}^{m-1} (p^*(0,k) + p^*(0,-k))$$

$$\geq \frac{4}{\sqrt{2\pi e}} \sum_{k=1}^{m-1} \frac{1}{k} - 2C \sum_{k=1}^{m-1} \frac{1}{k^3} \geq \frac{4}{\sqrt{2\pi e}} \log n - C'$$

for a universal constant $C' > 0$. On the other hand, the stationary probabilities are all $\frac{1}{n}$, so Proposition 1.5 gives the bound
$$\sum_{k \in \mathbb{Z}/n} p^*(0, k) \leq 2e \cdot \log n.$$
Thus, Proposition 1.5 is tight up to constant factor for random walk on the cycle.

**Remark 2.3.** It would be interesting to determine whether the aforementioned constant factor can be removed for random walks on regular graphs. More precisely, is it true that for any $\epsilon > 0$ and sufficiently large $n$, for a random walk on an $n$-vertex graph we have
$$\sum_{y \in S} p^*(x, y) \leq \left(\frac{4}{\sqrt{2\pi e}} + \epsilon\right) \log n \quad ?$$

## 2.2 Lower bound for Proposition 1.1

The following construction shows that Proposition 1.1 is optimal up to constant factor.

**Claim 2.1.** *For any $n \geq 2$ and $t \geq 2n$, there exists a Markov chain with state space $S$ of size $n$ and states $x, y \in S$ such that*
$$\mathbf{P}_x(\tau(y) = t) \geq \frac{n}{8t}.$$

*Proof.* The following construction is due to Kozma and Zeitouni (private communication). Write $t = r(n-1) + k$, where $1 \leq k \leq n-1$ and $r \geq 2$. Consider the Markov chain with states $s_1, s_2, \ldots, s_{n-1}$, and one additional state $u$ (see figure 1). We take the transition probabilities to be $p(s_i, s_{i+1}) = 1$ for each $1 \leq i \leq n-2$, and
$$p(s_{n-1}, s_1) = 1 - q, \quad p(s_{n-1}, u) = q, \quad p(u, s_2) = 1,$$
where $q = \frac{1}{r}$. Note that this chain is periodic with period $n - 1$.

We consider the hitting time from $s_{n-k}$ to $u$. Note that in the first $k$ steps the process deterministically goes to $s_{n-1}$ and then goes to either $u$ or $s_1$ with probabilities $q$ and $1-q$, respectively. Thereafter, every $n-1$ steps it again goes to either $u$ or $s_1$ with those probabilities. Thus,
$$\mathbf{P}_{s_{n-k}}(\tau(u) = t) = \left(1 - \frac{1}{r}\right)^r \cdot \frac{1}{r} \geq \frac{1}{4r} = \frac{n-1}{4t - 4k} \geq \frac{n}{8t},$$
completing the construction for $x = s_{n-k}$ and $y = u$. □

## 2.3 Lower bound for Theorem 1.2

The next construction is a slight modification of a construction made in [8], where it was given as an example of cutoff in the mixing time of trees. It gives the lower bound (1) mentioned in the introduction.

**Claim 2.2.** *There exist simple graphs of $n$ vertices, for arbitrarily large values of $n$, such that for the random walk started at a vertex $x$, there is another vertex $y$ and a time $t$ for which*
$$\mathbf{P}_x(\tau(y) = t) \geq \frac{c\sqrt{\log n}}{t},$$
*where $c > 0$ is a universal constant.*

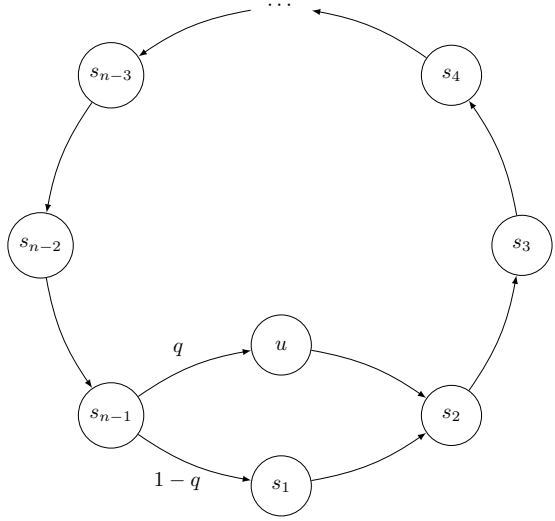

Figure 1: Illustration of Claim 2.1.

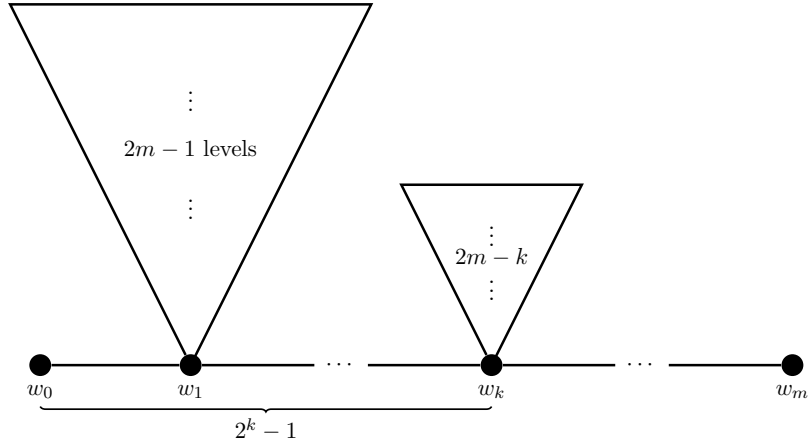

Figure 2: Illustration of $G_m$.

**Remark 2.4.** In fact, the graphs we construct have maximum degree 4 (corresponding to the vertices $w_k$ defined below). By adding self-loops, we can easily modify these examples to be regular graphs.

For any integer $m > 1$, consider the graph $G_m$ formed by a path of length $2^m - 1$ with vertices denoted by $v_1, v_2, \ldots, v_{2^m}$, together with $m - 1$ attached binary trees. In particular, for each $1 \leq k \leq m - 1$, we attach a binary tree of height $2m - k$ rooted at the vertex $w_k := v_{2^k}$ (see figure 2). The analysis of $\mathbf{P}_x(\tau(y) = t)$ for random walk on $G_m$ hinges on the following concentration of hitting times result similar to Lemma 2.1 of [8].

**Lemma 2.5.** *For the simple random walk on $G_m$ started at $w_m$, we have*

$$\mathbf{E}_{w_m}\tau(w_0) = \Theta\left(m \cdot 2^{2m}\right) \quad \text{and} \quad \mathbf{Var}_{w_m}(\tau(w_0)) = O\left(m \cdot 2^{4m}\right).$$

### 2.3.1 Proof of Lemma 2.5

First, we recall a standard estimate on trees.

**Lemma 2.6** ([8], Claim 2.3). *Let $Y_k$ denote the return time of a random walk on a binary tree of height $k$ started at the root. Then,*

$$\mathbf{E}Y_k = \Theta(2^k) \quad \text{and} \quad \mathbf{E}Y_k^2 = \Theta(2^{2k}).$$

Let $\mathcal{T}_k$ denote the tree (of height $2m - k$) attached to $w_k$ (excluding $w_k$ itself), and let $\mathcal{F} = \bigcup_{k=1}^{m} \mathcal{T}_k$. For $1 \leq k \leq m$, consider the walk started at $w_k$ and stopped upon hitting $w_{k-1}$. Define $Z_k$ to be the number of steps either starting or ending in $\mathcal{F}$. Note that $Z_m$ is deterministically zero. We have the following moment bounds on $Z_k$.

**Lemma 2.7.** *We have for $1 \leq k \leq m-1$,*

$$\mathbf{E}Z_k = \Theta(2^{2m}) \quad \text{and} \quad \mathbf{E}Z_k^2 = \Theta(2^{4m}).$$

*Proof.* The path of the random walk started at $w_k$ and stopped upon hitting $w_{k-1}$ may be decomposed into components of the following types:

(i) excursions from $w_k$ hitting $w_{k+1}$ and not hitting $\mathcal{T}_k$

(ii) excursions from $w_k$ staying entirely within $\mathcal{T}_k$

(iii) excursions from $w_k$ not hitting $\{w_{k-1}, w_{k+1}\} \cup \mathcal{T}_k$

(iv) a path from $w_k$ to $w_{k-1}$

Note that components of types (iii) and (iv) do not contribute to the time spent in $\mathcal{F}$. Components of type (i) can only spend time in $\mathcal{F}$ after first hitting $w_{k+1}$, so the time they spend in $\mathcal{F}$ is distributed as $Z_{k+1}$. The time spent in $\mathcal{F}$ by a component of type (ii) is like an excursion from the root for $\mathcal{T}_k$, so it is distributed as $Y_{2m-k}$. Thus, the law of $Z_k$ is given by the sum of (a random number of) independent copies of $Z_{k+1}$ and $Y_{2m-k}$.

To compute the distribution of how many components of types (i) and (ii) there are, let $\tau_k$ be the first time that the random walk started at $w_k$ hits $\{w_{k-1}, w_{k+1}\} \cup \mathcal{T}_k$. By a standard calculation,

$$\mathbf{P}_{w_k}(X_{\tau_k} = w_{k-1}) = \frac{2^{-k}}{2^{-k} + 2^{-k-1} + 1}$$

$$\mathbf{P}_{w_k}(X_{\tau_k} = w_{k+1}) = \frac{2^{-k-1}}{2^{-k} + 2^{-k-1} + 1}$$

$$\mathbf{P}_{w_k}(X_{\tau_k} \in \mathcal{T}_k) = \frac{1}{2^{-k} + 2^{-k-1} + 1}$$

It follows that the number $K$ of components of type (i) or (ii) is geometrically distributed with mean

$$\left(\frac{2^{-k}}{2^{-k} + 2^{-k-1} + 1}\right)^{-1} = 2^k + 1 + \frac{1}{2}.$$

Let $W_k$ be a random variable drawn according to $Z_{k+1}$ with probability $\frac{1}{2^{k+1}+1}$ and according to $Y_{2m-k}$ with probability $\frac{2^{k+1}}{2^{k+1}+1}$. Then, conditioned on a component being of type (i) or (ii), the time it spends in $T$ is distributed as $W_k$. Thus, $Z_k$ is a sum of $K$ i.i.d. copies of $W_k$. We now give recursive estimates of the first and second moments of $Z_k$. For the first moment, we have

$$\mathbf{E}Z_k = (\mathbf{E}K) \cdot (\mathbf{E}W_k) = (\mathbf{E}K) \cdot \left( \frac{2^{k+1}}{2^{k+1}+1} \mathbf{E}Y_{2m-k} + \frac{1}{2^{k+1}+1} \mathbf{E}Z_{k+1} \right).$$

Applying Lemma 2.6 to estimate $\mathbf{E}Y_{2m-k}$ and using our explicit characterization of $K$, this gives

$$c \cdot 2^{2m} \leq \mathbf{E}Z_k \leq C \cdot 2^{2m} + \frac{7}{10} \cdot \mathbf{E}Z_{k+1}$$

for some universal constants $c$ and $C$. Recall that $Z_m = 0$, so this recursive bound implies that $\mathbf{E}Z_k = \Theta(2^{2m})$ for all $1 \leq k \leq m-1$.

We now bound the second moment of $Z_k$. Note that our first moment bound immediately implies that $\mathbf{E}W_k = \Theta(2^{2m-k})$. We have

$$\mathbf{E}Z_k^2 = (\mathbf{E}K) \cdot (\mathbf{E}W_k^2) + (\mathbf{E}K(K-1)) \cdot (\mathbf{E}W_k)^2$$

$$= (\mathbf{E}K) \cdot \left( \frac{2^{k+1}}{2^{k+1}+1} \mathbf{E}Y_{2m-k}^2 + \frac{1}{2^{k+1}+1} \mathbf{E}Z_{k+1}^2 \right) + O(2^{4m})$$

$$\leq \frac{7}{10} \cdot \mathbf{E}Z_{k+1}^2 + O(2^{4m}),$$

where we use the upper bound on $\mathbf{E}Y_{2m-k}^2$ from Lemma 2.6. By the same argument as was used for the first moment, this implies $\mathbf{E}Z_k^2 = O(2^{4m})$ for $1 \leq k \leq m-1$. Note that we also have $\mathbf{E}Z_k^2 \geq (\mathbf{E}Z_k)^2$, so we have

$$\mathbf{E}Z_k = \Theta(2^{2m}) \quad \text{and} \quad \mathbf{E}Z_k^2 = \Theta(2^{4m}).$$

as desired. $\square$

We are now ready to prove Lemma 2.5.

*Proof of Lemma 2.5.* Note that a random walk started at $w_m$ and stopped upon hitting $w_0$ can be decomposed into independent walks started at $w_k$ and stopped upon hitting $w_{k-1}$ for $k = m, m-1, \ldots, 1$. Let $Z = Z_1 + Z_2 + \cdots + Z_m$ be the total number of steps starting or ending in $\mathcal{F}$, and let $Z'$ be the total number of steps completely outside of $\mathcal{F}$, so that $\tau := Z + Z'$ is the hitting time from $w_m$ to $w_0$.

Note that $Z'$ has the distribution of the hitting time of a random walk on a path of length $2^m - 1$ from one end to the other. The following estimates are standard (see e.g. [2], III.7):

$$\mathbf{E}Z' = \Theta(2^{2m}) \quad \text{and} \quad \mathbf{E}Z'^2 = O(2^{4m}).$$

It follows from Lemma 2.7 that

$$\mathbf{E}\tau = \mathbf{E}Z + \mathbf{E}Z' = \Theta(m \cdot 2^{2m}),$$

which proves the first part of the lemma. We also have

$$\mathbf{Var}(Z) = \sum_{k=1}^{m} \mathbf{Var}(Z_k) = O(m \cdot 2^{4m})$$

$$\mathbf{E}Z^2 = (\mathbf{E}Z)^2 + \mathbf{Var}(Z) = O(m^2 \cdot 2^{4m}).$$

Thus,
$$\mathbf{Var}(\tau) = 2\left(\mathbf{E}ZZ' - \mathbf{E}Z\mathbf{E}Z'\right) + \mathbf{Var}(Z) + \mathbf{Var}(Z')$$
$$\leq 2\sqrt{(\mathbf{E}Z^2)(\mathbf{E}Z'^2)} + \mathbf{Var}(Z) + \mathbf{Var}(Z') = O(m \cdot 2^{4m}),$$

which is the second statement in the lemma, completing the proof. □

### 2.3.2 Proof of Claim 2.2

We will show using Lemma 2.5 that $G_m$ satisfies the criteria of Claim 2.2.

*Proof of Claim 2.2.* For any $m \geq 2$, take
$$n = 2^{2m} - 2^m - 2m + 2$$
to be the number of vertices of $G_m$. For convenience, write $T = \mathbf{E}_{w_m}\tau(w_0)$. By Lemma 2.5 and Chebyshev's inequality, there is a universal constant $C$ such that
$$\mathbf{P}_{w_m}\left(|\tau(w_0) - T| \geq C\sqrt{m \cdot 2^{4m}}\right) \leq \frac{1}{2}.$$

It follows by the pigeonhole principle that for some $t$ with
$$T - C\sqrt{m \cdot 2^{4m}} \leq t \leq T + C\sqrt{m \cdot 2^{4m}} \tag{3}$$
we have
$$\mathbf{P}_{w_m}\left(\tau(w_0) = t\right) \geq \frac{1}{4C\sqrt{m} \cdot 2^{2m}}.$$

Moreover, note that $T$ is on the order of $m \cdot 2^{2m}$ by Lemma 2.5, so for sufficiently large $m$, (3) implies that $t \geq cm \cdot 2^{2m}$ for a sufficiently small constant $c > 0$. Then,
$$\mathbf{P}_{w_m}\left(\tau(w_0) = t\right) \geq \frac{1}{4C\sqrt{m} \cdot 2^{2m}} \geq \frac{c\sqrt{m}}{4Ct} \geq \frac{c}{4\sqrt{2} \cdot C} \cdot \frac{\sqrt{\log n}}{t}.$$

This is the desired bound upon renaming of constants, and $m$ (hence $n$) can be made arbitrarily large. □

## 3 Sums of geometrics

We devote this section to establishing some basic estimates for geometric random variables that will be used later. This section can be read on its own and does not rely on any other part of the paper. We say that a random variable $Z$ is *geometrically distributed with parameter $p$* if $\mathbf{P}(Z = t) = p(1-p)^t$ for $t = 0, 1, 2, \ldots$. First, we have the following estimate on sums of i.i.d. geometrics.

**Lemma 3.1.** *Let $m > 0$ be an integer, and let $X_1, X_2, \ldots, X_n$ be i.i.d. geometric random variables with parameter $p = \frac{n}{m+n}$. Then,*
$$\frac{1}{3}\sqrt{\frac{n}{m(m+n)}} \leq \mathbf{P}\left(X_1 + X_2 + \cdots + X_n = m\right) \leq \frac{1}{2}\sqrt{\frac{n}{m(m+n)}}.$$

For the proof, we use the following version of Stirling's formula which holds for all integers $N \geq 1$ (see e.g. [2], II.9).

$$\sqrt{2\pi} \leq \frac{N!}{N^{N+\frac{1}{2}}e^{-N}} \leq \sqrt{2\pi} \cdot e^{\frac{1}{12}}. \tag{4}$$

Consequently, we have the following estimate on binomial coefficients.

**Proposition 3.2.** *For any positive integers $M$ and $N$, we have*

$$\frac{1}{3}\sqrt{\frac{M+N}{MN}} \cdot \frac{(M+N)^{M+N}}{M^M \cdot N^N} \leq \binom{M+N}{M} \leq \frac{1}{2}\sqrt{\frac{M+N}{MN}} \cdot \frac{(M+N)^{M+N}}{M^M \cdot N^N}.$$

*Proof.* Apply (4) to each factorial term in $\binom{M+N}{M} = \frac{(M+N)!}{M! \cdot N!}$. □

*Proof of Lemma 3.1.* By direct calculation, we have

$$\mathbf{P}(X_1 + X_2 + \cdots + X_n = m) = p^n(1-p)^m \binom{m+n-1}{n-1}$$

$$= \left(\frac{n}{m+n}\right)^n \left(\frac{m}{m+n}\right)^m \cdot \binom{m+n}{n} \cdot \frac{n}{m+n}.$$

Applying Proposition 3.2 to $\binom{m+n}{n}$ gives the result. □

We can also prove the following bound on sums of independent geometric variables (not necessarily identically distributed) due to Nazarov (private communication).

**Lemma 3.3.** *Let $X_1, X_2, \ldots, X_n$ be independent geometric random variables. Then, for any $t \geq 0$,*

$$\mathbf{P}(X_1 + X_2 + \cdots + X_n = t) \leq \frac{1}{2}\sqrt{\frac{n}{t(t+n)}}. \tag{5}$$

**Remark 3.4.** We actually show that $\mathbf{P}(X_1 + X_2 + \cdots + X_n = t)$ is maximized when each $X_i$ has parameter $\frac{n}{t+n}$. The lemma then follows from Lemma 3.1.

*Proof.* Let $p_i$ be the parameter of $X_i$, and define $q_i = 1 - p_i \in [0, 1)$, so that $\mathbf{P}(X_i = k) = (1-q_i)q_i^k$. Then, by expanding terms, we have

$$\mathbf{P}(X_1 + X_2 + \cdots + X_n = t) = \prod_{i=1}^{n}(1-q_i) \sum_{\substack{j_1+\ldots+j_n=t \\ j_i \in \mathbb{Z}^+}} q_1^{j_1} q_2^{j_2} \cdots q_n^{j_n}. \tag{6}$$

The right side of equation (6) may be regarded as a function of $q = (q_1, \ldots, q_n) \in [0,1]^n$, which we denote by $F(q)$.

We will show that $F(q)$ is maximized when all the $q_i$ are equal. This will be accomplished by showing that if $q_i \neq q_j$ for some $i$ and $j$, then by replacing both $q_i$ and $q_j$ with some common value, we can increase $F(q)$. To this end, for any $0 \leq x < y < 1$ and integer $k \geq 0$, define

$$f_k(x,y) = (1-x)(1-y)(x^k + x^{k-1}y + \ldots + y^k).$$

We also define

$$P(x,y,t) = \begin{cases} \frac{(1-x)(1-y)}{(y-x)(1-t)^2} & \text{if } x \leq t \leq y \\ 0 & \text{otherwise.} \end{cases}$$

Clearly, $P(x, y, t) \geq 0$, and

$$\int_0^1 P(x, y, t)\, dt = \int_x^y P(x, y, t)\, dt = \frac{(1-x)(1-y)}{(y-x)(1-t)}\bigg|_{t=x}^y = 1,$$

and so $P(x, y, \cdot)$ may be regarded as a probability density on $[0, 1]$. Moreover,

$$\int_0^1 P(x, y, t) f_k(t, t)\, dt = \int_x^y \frac{(1-x)(1-y)}{y-x} \cdot (k+1) t^k\, dt$$

$$= \frac{(1-x)(1-y)}{y-x} \cdot (y^{k+1} - x^{k+1}) = f_k(x, y). \tag{7}$$

Now, consider a point $q^* = (q_1^*, \ldots, q_n^*) \in [0,1]^n$ which maximizes $F(q^*)$ on the compact set $[0,1]^n$, and suppose for sake of contradiction that not all of the $q_i^*$ are equal. Then, without loss of generality, we may assume $q_1^* < q_2^*$. Fixing $q_i^*$ for $i > 2$ and considering $F$ as a function of the first two variables only, it is easy to see from equation (6) that $F$ takes the form

$$F(q_1^*, q_2^*) = \sum_{k=0}^t c_k f_k(q_1^*, q_2^*)$$

for constants $c_0, c_1, \ldots, c_t$. Thus, according to (7) with $x = q_1^*$ and $y = q_2^*$, we have[1]

$$F(q_1^*, q_2^*) = \int_0^1 P(q_1^*, q_2^*, t) F(t, t)\, dt.$$

However, by the maximality of $F(q^*)$, we also have that $F(t, t) \leq F(q_1^*, q_2^*)$ for all $t \in [q_1^*, q_2^*]$, so

$$\int_0^1 P(q_1^*, q_2^*, t) F(t, t)\, dt \leq \int_0^1 P(q_1^*, q_2^*, t) F(q_1^*, q_2^*)\, dt = F(q_1^*, q_2^*),$$

and this inequality is strict since $F(t, t)$ is non-constant on $[q_1^*, q_2^*]$. This is a contradiction, so it follows that the $q_i^*$ all take some common value $r$. To determine this value $r$, we may compute

$$F(q^*) = (1-r)^n \cdot \binom{t+n-1}{n-1} r^t,$$

and optimizing over $r$, we find that $r = \frac{t}{t+n}$. This corresponds to the case where each $X_i$ has parameter $\frac{n}{t+n}$. Applying Lemma 3.1 finishes the proof. □

**Corollary 3.5.** *If $(X_i)_{i=1}^n$ are independent mixtures of geometric random variables, then (5) holds.*

*Proof.* For each $i$, let $\theta_i$ be a random variable on $[0, 1]$ so that $X_i$ is distributed as a geometric with parameter $\theta_i$. By Lemma 3.3, we have

$$\mathbf{P}(X_1 + X_2 + \cdots + X_n = t \mid \theta_1, \ldots, \theta_n) \leq \frac{1}{2}\sqrt{\frac{n}{t(t+n)}}.$$

Taking the expectation over the $\theta_i$ gives the result. □

---

[1]The assumption $y < 1$ is satisfied because clearly if $q_2^* = 1$, then $F(q^*) = 0$ is not maximal.

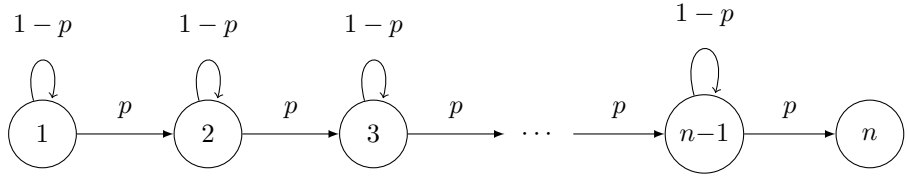

Figure 3: Illustration of Claim 4.1.

## 4 Proof of Theorem 1.3

We are now ready to prove Theorem 1.3, which is a "worst case" bound for reversible chains that does not depend on the stationary probability. Before proving the theorem, it is instructive to construct the example that attains the lower bound.

**Claim 4.1.** *For any $n \geq 2$ and $t \geq n$, there exists a reversible Markov chain $X$ on $n$ states with two states $x$ and $y$ such that*
$$\mathbf{P}_x(\tau(y) = t) \geq \frac{\sqrt{n}}{3t}.$$

*Proof.* Let $p = \frac{n}{t}$. Consider the case where $X$ is a pure-birth chain with states labeled $1, 2, \ldots, n$, with transition probabilities
$$p(i, i) = 1 - p, \quad p(i, i+1) = p \text{ for } 1 \leq i \leq n-1$$
$$p(n, n) = 1.$$

Suppose that $X$ starts at 1, and let $T_i$ be the first hitting time of state $i$, with $T_1 = 0$. Let $D_i = T_{i+1} - T_i - 1$, so that
$$T_n = n - 1 + \sum_{i=1}^{n-1} D_i.$$

Note that each $D_i$ is a geometric random variable with parameter $p$, and the $D_i$ are independent. Thus, applying Lemma 3.1 with $m = t - n + 1$, we have
$$\mathbf{P}(T_n = t) \geq \frac{1}{3}\sqrt{\frac{n}{(t-n+1)(t+1)}} \geq \frac{\sqrt{n}}{3t},$$

as desired. □

In fact, the above example is in some sense the extremal case. In a pure-birth chain, we saw that the hitting time from the beginning to the end is a sum of geometric random variables. Lemma 3.3 implies that the sum of $n$ geometric random variables cannot be too concentrated, which essentially proves Theorem 1.3 for pure-birth chains.

Then, we will show that the behavior of the hitting time in a general reversible Markov chain is like a mixture of pure-birth chains. Such a representation was shown by Miclo [7]; in this paper, we give a short elementary proof based on loop erasures of first hitting paths from $x$ to $y$. By this mixture argument, it follows that the bound for pure-birth chains carries over to the general case, although we lose approximately a factor of 2 due to the possibility of negative eigenvalues.

The next lemma is a variant of the well-known spectral decomposition of return probabilities in reversible Markov chains. The proof is essentially the same as that of Lemma 2.2 in [3], although our formulation includes an additional non-negativity statement. For the sake of completeness, we include a proof of this lemma in Section 9.2.

**Lemma 4.1.** *Let $X$ be a reversible, irreducible Markov chain with finite state space $S$. Consider any $x \in S$, and let $U \subset S$ be a subset not containing $x$. Then, there exist real numbers $a_1, \ldots, a_{|S|} \geq 0$ and $\lambda_1, \ldots, \lambda_{|S|} \in [-1, 1]$ such that*

$$\mathbf{P}_x(X_t = x, \tau(U) > t) = \sum_{i=1}^{|S|} a_i \lambda_i^t.$$

*for each integer $t \geq 0$. Moreover, if the eigenvalues of $X$ (that is, the eigenvalues of the matrix of transition probabilities) are non-negative, then we may take $\lambda_i \in [0, 1]$ for each $i$.*

Let us now prove Theorem 1.3.

*Proof of Theorem 1.3.* We will first prove the stronger bound

$$\mathbf{P}_x(\tau(y) = t) \leq \frac{1}{2}\sqrt{\frac{n}{t(t-n)}}. \tag{8}$$

in the case where all the eigenvalues are non-negative. We then deduce the theorem by considering even times. It is convenient to assume throughout that $X$ is irreducible so Lemma 4.1 applies. This is valid because we may always restrict to the communicating class of $x$ without increasing the number of states.

Let $z = (x = z_0, z_1, \ldots, z_k = y)$ be a random variable denoting the path taken by the chain started at $x$ and stopped upon hitting $y$. Thus, $|z| = \tau_x(y)$, and we are interested in bounding $\mathbf{P}(|z| = t)$. Define the *loop erasure* of $z$ to be the path $(w_0, w_1, \ldots, w_\ell)$ determined as follows. We take $w_0 = x$, and inductively for each $i \geq 0$, let $k_i$ be the largest index such that $z_{k_i} \in \{w_0, \ldots, w_i\}$. Then, as long as $k_i < k$, define $w_{i+1} = z_{k_i+1}$ and continue the process for $i + 1$. If $k_i = k$, then the path ends. In particular, $k_\ell = k$.

We denote this path, which is a function of $z$, by $w(z)$. A less formal description of the loop erasure is that $w(z)$ is the path obtained by following $z$ and, whenever the path forms a loop, removing all vertices from that loop. Loop-erased walks appear in other contexts, including random walks on lattices and uniform sampling of spanning trees (see e.g. [4], [10]).

Fix a loop erasure $w = (x = w_0, w_1, \ldots, w_\ell = y)$, and let $\hat{\mathbf{P}}_w$ denote the conditional probability $\mathbf{P}(\cdot \mid w(z) = w)$. We now describe a method of sampling from $\hat{\mathbf{P}}_w(z)$. For each $0 \leq i < \ell$, let $\mathcal{P}_i$ denote the set of all paths starting and ending at $w_i$ (possibly of length 1) and avoiding $\{y, w_0, w_1, \ldots, w_{i-1}\}$. For each $i$, we independently sample a path $\tilde{z}_i = (\tilde{z}_{i,0}, \tilde{z}_{i,1}, \ldots \tilde{z}_{i,k})$ from $\mathcal{P}_i$ with probability

$$c_i \cdot \prod_{j=1}^{k} p(\tilde{z}_{i,j-1}, \tilde{z}_{i,j}),$$

where $c_i$ is the normalizing constant which makes all the probabilities sum to 1. We obtain a sample from $\hat{\mathbf{P}}_w(z)$ by taking the concatenation $z = [\tilde{z}_0][\tilde{z}_1] \cdots [\tilde{z}_{\ell-1}]y$.

Now, observe that

$$\mathbf{P}(|\tilde{z}_i| = t) = \sum_{\substack{w \in \mathcal{P}_i \\ |w|=t}} c_i \cdot \prod_{j=1}^{t} p(w_{j-1}, w_j) = c_i \cdot \mathbf{P}_{w_i}(X_t = w_i, \tau(\{y, w_1, \ldots, w_{i-1}\}) > t).$$

By Lemma 4.1, $|\tilde{z}_i|$ is distributed as a mixture of geometric random variables. Also, note that

$$|z| = \ell + |\tilde{z}_0| + |\tilde{z}_1| + \cdots + |\tilde{z}_{\ell-1}|.$$

Hence, by Corollary 3.5, we have

$$\hat{\mathbf{P}}_w(|z| - \ell = m) \leq \frac{1}{2}\sqrt{\frac{\ell}{m(m+\ell)}}$$

Recall that $w$ is a loop erasure and therefore, the $w_i$ are all distinct, which implies $\ell \leq n$. Thus, taking $m = t - \ell$, the above inequality can be rewritten

$$\hat{\mathbf{P}}_w(|z| = t) \leq \frac{1}{2}\sqrt{\frac{\ell}{(t-\ell)t}} \leq \frac{1}{2}\sqrt{\frac{n}{(t-n)t}}.$$

This bound holds for each $w$, so taking the expectation over all possible loop erasures $w$, we obtain

$$\mathbf{P}(|z| = t) \leq \frac{1}{2}\sqrt{\frac{n}{(t-n)t}}.$$

We have thus proved (8) when the chain has non-negative eigenvalues. Let us now consider the general case, where the eigenvalues can be negative. Let $Y_t = X_{2t}$, so that $Y$ is also a Markov chain (with transition matrix the square of the transition matrix of $X$). Note that $Y$ has non-negative eigenvalues. Let $\tau_{x,X}(y)$ denote the hitting time from $x$ to $y$ under the chain $X$, and similarly let $\tau_{x,Y}(y)$ be the hitting time from $x$ to $y$ under $Y$. If $\tau_{x,X}(y) = 2k$, then we necessarily have $\tau_{x,Y}(y) = k$. Hence, if $t = 2k$, we immediately have

$$\mathbf{P}(\tau_{x,X}(y) = t) \leq \mathbf{P}(\tau_{x,Y}(y) = k) \leq \frac{1}{2}\sqrt{\frac{n}{(k-n)k}} < \frac{\sqrt{2n}}{t},$$

where we used the assumption $t \geq 4n + 4$. If $t = 2k + 1$, then we have

$$\mathbf{P}(\tau_{x,X}(y) = t) = \sum_{s \in S,\, s \neq y} p(x,s) \mathbf{P}(\tau_{s,X}(y) = t-1)$$

$$\leq \sum_{s \in S,\, s \neq y} p(x,s) \cdot \frac{1}{2}\sqrt{\frac{n}{(k-n)k}} \leq \frac{1}{2}\sqrt{\frac{n}{(k-n)k}} \leq \frac{\sqrt{2n}}{t}.$$

completing the proof. $\square$

## 5 Lower bound constructions for $\mathbf{P}_x(\mathcal{S}_t)$

We now describe how to modify the constructions in Claims 2.1, 2.2, and 4.1 to give corresponding lower bounds for $\mathbf{P}_x(\mathcal{S}_t)$.

First, recall the cycle construction in the proof of Claim 2.1. Using a similar construction but with $n$ copies of the state $u$ (see figure 4), we can obtain the following lower bound for $\mathbf{P}_x(\mathcal{S}_t)$.

**Claim 5.1.** *For any $n \geq 2$ and $t \geq 2n$, there exists a Markov chain with $2n$ states and a starting state $x$ such that*

$$\mathbf{P}_x(\mathcal{S}_t) \geq \frac{n^2}{56t}.$$

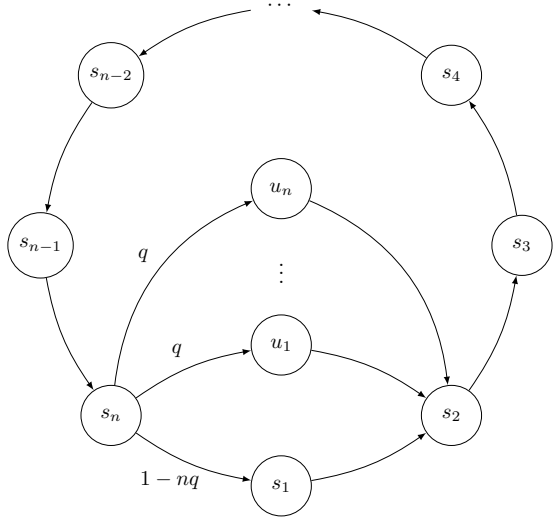

Figure 4: Illustration of Claim 5.1.

*Proof.* Write $t = rn + k$, where $1 \leq k \leq n$ and $r \geq 2$. Define $q = \frac{1}{r}$. Consider a Markov chain with states
$$\{s_1, s_2, \ldots, s_n\} \cup \{u_1, u_2, \ldots, u_n\}$$
whose transition probabilities are

$$p(s_i, s_{i+1}) = 1 \quad \text{for each } 1 \leq i \leq n-1,$$
$$p(s_n, s_1) = 1 - nq$$
$$p(s_n, u_i) = q, \quad p(u_i, s_2) = 1 \quad \text{for each } 1 \leq i \leq n.$$

This is very similar to the example in the proof of Claim 2.1, except that $u$ is replaced by many states $u_1, u_2, \ldots, u_n$. Note that the chain is periodic with period $n$.

We consider the process started at $x := s_{n-k+1}$, so that at time $t - 1$, the process is in state $s_n$. After $t - 1$ steps, the process has taken $r$ steps starting in $s_n$, so the probability that a given state $u_j$ is still unvisited is $(1 - q)^r$. Thus, letting $Z$ be the number of states in $\{u_1, \ldots, u_n\}$ which are visited by time $t - 1$, we have

$$\mathbf{E}Z = n \cdot (1 - (1-q)^r) = n \cdot \left(1 - \left(1 - \frac{1}{r}\right)^r\right) \leq \frac{3}{4} \cdot n,$$

and so Markov's inequality implies

$$\mathbf{P}_x\left(Z \geq \frac{7n}{8}\right) \leq \frac{6}{7}.$$

On the event that $Z < \frac{7n}{8}$, there are at least $\frac{n}{8}$ states among $\{u_1, \ldots, u_n\}$ which have not been visited, so the probability that the chain hits a new state on the $t$-th step is at least $\frac{nq}{8}$. Thus,

$$\mathbf{P}_x(\mathcal{S}_t) \geq \mathbf{P}_x\left(\mathcal{S}_t \,\bigg|\, Z < \frac{7n}{8}\right) \cdot \mathbf{P}_x\left(Z < \frac{7n}{8}\right) \geq \frac{nq}{56} \geq \frac{n^2}{56t},$$

as desired. □

We can also modify the constructions corresponding to Claims 2.2 and 4.1 to obtain the following.

**Claim 5.2.** *There exist simple graphs of n vertices, for arbitrarily large values of n, such that for the random walk started at a vertex x, there is a time t for which*

$$\mathbf{P}_x(\mathcal{S}_t) \geq \frac{cn\sqrt{\log n}}{t},$$

*where $c > 0$ is a universal constant.*

**Claim 5.3.** *For any n, there exists a reversible Markov chain X on 2n states with two states x and y such that*

$$\mathbf{P}_x(\mathcal{S}_t) = \Omega\left(\frac{n\sqrt{n}}{t}\right).$$

The constructions we give for the above claims are both based on a general lemma that translates certain lower bounds on $\mathbf{P}_x(\tau(y) = t)$ into lower bounds on $\mathbf{P}_x(\mathcal{S}_t)$.

**Lemma 5.1.** *Consider a Markov chain with state space S with states $x, y \in S$ and a subset $U \subset S$ such that starting from x, the chain cannot reach U without going through y. Let $Z_s$ denote the number of visited states in U at time s. Then, for any integer $N > 0$, there exists s with $t \leq s < t + 2N$ such that*

$$\mathbf{P}_x(\mathcal{S}_s) \geq \frac{\mathbf{P}_x(t \leq \tau(y) < t + N) \cdot \mathbf{E}_y Z_N}{2N}.$$

*Proof.* Note that no state in $U$ can be visited before $\tau(y)$. We lower bound $\mathbf{E}_x(Z_{t+2N-1} - Z_t)$ by only considering the event that $t \leq \tau(y) < t + N$. We have

$$\mathbf{E}_x(Z_{t+2N-1} - Z_t) \geq \sum_{k=0}^{N-1} \mathbf{P}_x(\tau(y) = t+k) \mathbf{E}_y Z_N = \mathbf{P}_x(t \leq \tau(y) < t+N) \cdot \mathbf{E}_y Z_N.$$

This lower bounds the expected number of new states visited in the time interval $[t, t+2N)$, so it follows by the pigeonhole principle that for some $t \leq s < t + 2N$, we have

$$\mathbf{P}_x(\mathcal{S}_s) \geq \frac{\mathbf{P}_x(t \leq \tau(y) < t + N) \cdot \mathbf{E}_y Z_N}{2N}.$$

□

We now prove the claims.

*Proof Claim 5.2.* Let $G_m$ be as in the proof of Claim 2.2, and let $n$ be the size of $G_m$. Let $k = \lfloor \sqrt[3]{n} \rfloor$, and let $H_n$ denote the three-dimensional discrete torus of size $k^3$, whose vertex set is $(\mathbb{Z}/k\mathbb{Z})^3$, with edges between nearest neighbors.

We recall the standard fact that the effective resistance between $y$ and $z$ in $H_k$ is bounded by a constant (Exercise 9.1, [6]). Define a graph $\tilde{G}_m$ by attaching to $G_m$ a copy of $H_k$ so that $(0,0,0)$ in $H_k$ is joined to $w_0$ in $G_m$ (see figure 5). The graph $\tilde{G}_m$ has $O(k^3) = O(n)$ edges, so by the commute time identity (Proposition 10.6, [6]), we have for any $z \in H_k$ (with $H_k$ regarded as a subgraph of $\tilde{G}_m$),

$$\mathbf{E}_{w_0} \tau(z) = O(n).$$

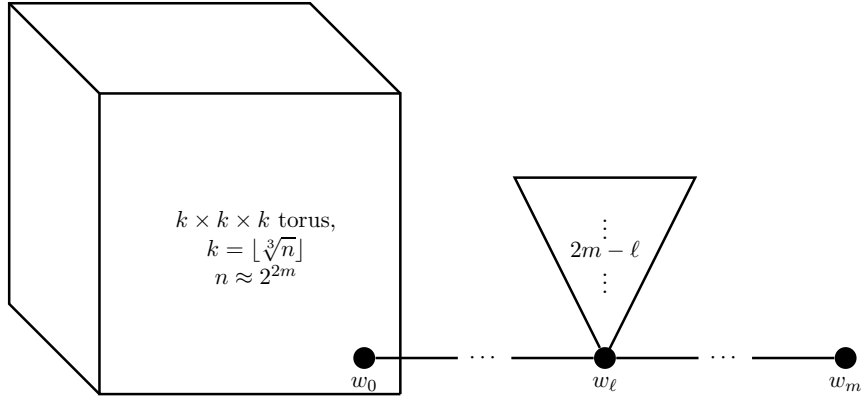

Figure 5: Illustration of $\tilde{G}_m$.

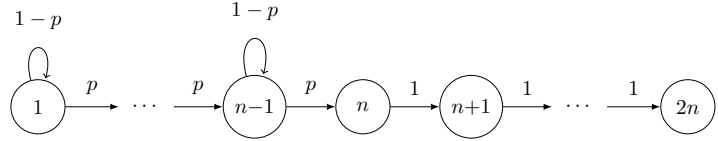

Figure 6: Illustration of Claim 5.3.

Consequently, the expected number of visited vertices in $H_n$ for a random walk started at $w_0$ and run for $n$ steps is at least $cn$ for some constant $c > 0$. Recall also from Lemma 2.5 that $\mathbf{E}_{w_m}\tau(w_0)$ has expectation on the order of $m \cdot 2^{2m} = \Theta(n \log n)$ with fluctuations on the order of $\sqrt{m} \cdot 2^{2m} = \Theta(n\sqrt{\log n})$. It follows that for some $t$ with $t = \Theta(n \log n)$,

$$\mathbf{P}_{w_m}(t \leq \tau(w_0) < t + n) = \Omega\left(\frac{1}{\sqrt{\log n}}\right).$$

Thus, we may apply Lemma 5.1 with $U$ being the vertex set of $H_k$, $x = w_m$, $y = w_0$, and $N = n$. We find that for some $s = \Theta(n\sqrt{\log n})$,

$$\mathbf{P}_{w_m}(\mathcal{S}_s) \geq \frac{\mathbf{P}_{w_m}(t \leq \tau(w_0) < t + n) \cdot \mathbf{E}_{w_0} Z_n}{2n}$$

$$= \Omega\left(\frac{1}{\sqrt{\log n}}\right) \cdot \Theta(n) \cdot \frac{1}{2n} = \Omega\left(\frac{n\sqrt{\log n}}{s}\right),$$

as desired. □

*Proof Claim 5.3.* It suffices to prove the bound for $t \geq n\sqrt{n}$, so we assume this in what follows.

Consider the pure-birth chain from the proof of Claim 4.1 (with states labeled $1, 2, \ldots, n$), and introduce $n$ additional states labeled $n + 1, \ldots, 2n$. We use the same transition probabilities as in the proof of Claim 4.1, with the modification that $p(i, i+1) = 1$ for $i = n, n+1, \ldots, 2n-1$ (see figure 6).

Recall from the proof of Claim 4.1 that the hitting time $\tau$ from 1 to $n$ is distributed as

$$n - 1 + \sum_{i=1}^{n-1} D_i,$$

where the $D_i$ are independent geometrics of parameter $p = \frac{n}{t}$. We can thus calculate

$$\mathbf{E}\tau = \Theta(t) \quad \text{and} \quad \mathbf{Var}(\tau) = O\left(\frac{t^2}{n}\right).$$

It follows that for some $t' = \Theta(t)$,

$$\mathbf{P}(t' \leq \tau < t' + n) = \Omega\left(\frac{n\sqrt{n}}{t}\right).$$

We can apply Lemma 5.1 with $U = \{n, n+1, \ldots, 2n\}$, $x = 1$, $y = n$, and $N = n$. Note that if the chain is started at $y$, then $Z_n = n$ deterministically in this case. Thus, we obtain for some $s = \Theta(t') = \Theta(t)$,

$$\mathbf{P}_x(\mathcal{S}_s) \geq \frac{\mathbf{P}_x(t' \leq \tau(y) < t' + n) \cdot \mathbf{E}_y Z_n}{2n}$$

$$= \Omega\left(\frac{n\sqrt{n}}{t}\right) \cdot n \cdot \frac{1}{2n} = \Omega\left(\frac{n\sqrt{n}}{s}\right),$$

as desired. $\square$

## 6 Proofs of Propositions 1.7 and 1.8

*Proof of Proposition 1.7.* Let $S$ be the state space of $X$. Note that for $t > 1$,

$$\mathbf{P}_\pi(\tau(y) = t) = \sum_{x \in S \setminus \{y\}} \mathbf{P}_x(\tau(y) = t-1) \mathbf{P}_\pi(X_1 = x)$$

$$\leq \sum_{x \in S} \pi(x) \mathbf{P}_x(\tau(y) = t-1) = \mathbf{P}_\pi(\tau(y) = t-1).$$

Thus, $\mathbf{P}_\pi(\tau(y) = t)$ is non-increasing in $t$, and so

$$\mathbf{P}_\pi(\tau(y) = t) \leq \frac{1}{t}\sum_{s=0}^{t-1} \mathbf{P}_\pi(\tau(y) = s) \leq \frac{1}{t}.$$

For any particular state $x \in S$, we have

$$\mathbf{P}_x(\tau(y) = t) \leq \frac{1}{\pi(x)}\mathbf{P}_\pi(\tau(y) = t) \leq \frac{1}{t\pi(x)}.$$

$\square$

*Proof of Proposition 1.8.* We have

$$\mathbf{P}_x(\tau(y) = t) \leq \sum_{z \in S} p^s(x,z) \mathbf{P}_z(\tau(y) = t - s)$$

$$\leq \sum_{z \in S,\, p^s(x,z) \geq \pi(z)} (p^s(x,z) - \pi(z)) \mathbf{P}_z(\tau(y) = t-s) + \sum_{z \in S} \pi(z) \mathbf{P}_z(\tau(y) = t-s)$$

$$\leq \sum_{z \in S,\, p^s(x,z) \geq \pi(z)} (p^s(x,z) - \pi(z)) \psi(t-s) + \sum_{z \in S} \pi(z) \frac{1}{t-s}$$

$$\leq d(s)\psi(t-s) + \frac{1}{t-s},$$

where between the second and third lines we used Proposition 1.7.

Suppose that $t > 2t_{\mathrm{mix}}(1/4)\lceil \log_2 n\rceil$, and take $s = \lfloor \frac{t}{2} \rfloor$. By Proposition 1.1, we may take $\psi(u) = \frac{n}{u}$. Note that $s \geq t_{\mathrm{mix}}(1/4)\lceil \log_2 n\rceil$, so

$$d(s) \leq \left(\frac{1}{2}\right)^{\lceil \log_2 n\rceil} \leq \frac{1}{n}.$$

Thus,

$$\mathbf{P}_x(\tau(y) = t) \leq \frac{1}{n} \cdot \frac{n}{t-s} + \frac{1}{t-s} = \frac{2}{t-s} \leq \frac{4}{t}.$$

$\square$

# 7 Open problems

We pose two open problems arising from our work. First, it is natural to try to close the gap between the bound in Theorem 1.2 and the corresponding example given in Claim 2.2. We suspect that Theorem 1.2 is not optimal and make the following conjecture.

**Conjecture 7.1.** *Let $X$ be a reversible Markov chain with finite state space, and let $\pi$ be the stationary distribution of $X$. Consider any two states $x$ and $y$. For all $t > 0$,*

$$\mathbf{P}_x(\tau(y) = t) = O\left(\frac{1}{t}\sqrt{\max\left(1, \log \frac{1}{\pi(x)}\right)}\right).$$

The second question relates to hitting times in an interval and is due to Holroyd. The example in Claim 2.1 shows that Proposition 1.1 is optimal up to a constant, but it has the very special property of being periodic with a large period. Consequently, a time $t$ for which $\mathbf{P}_x(\tau(y) = t)$ is large is surrounded by many times $t'$ near $t$ where $\mathbf{P}_x(\tau(y) = t') = 0$. This motivates the following conjecture of Holroyd.

**Conjecture 7.2.** *Let $X$ be a Markov chain with finite state space $S$ with $|S| = n$, and let $x$ and $y$ be any two states. Then, for all $t > n$,*

$$\mathbf{P}_x(t \leq \tau(y) \leq t + n) = O\left(\frac{n}{t}\right).$$

# 8 Acknowledgements

We are grateful to Aryeh Kontorovich for posing the problem that led to this work and to Fedor Nazarov for sharing with us his elegant proof of Lemma 3.3. We also thank Jian Ding, Ronen Eldan, Jonathan Hermon, Alexander Holroyd, Gady Kozma, Shayan Oveis Gharan, Perla Sousi, Peter Winkler, and Ofer Zeitouni for helpful discussions. Finally, we are grateful to Microsoft Research for hosting visits by James Norris and Alex Zhai, which made this work possible.

# 9 Background

## 9.1 Proof of Theorem 2.2

The following distillation of Starr's proof came from helpful discussions with Jonathan Hermon.

*Proof.* We consider the chain where the initial state $X_0$ is drawn from the stationary measure $\pi$. It suffices to prove the case where $f \in L^p(S, \pi)$ is non-negative, since otherwise we may replace $f$ by $|f|$ (noting that $|P^{2n}f| \leq P^{2n}|f|$). Thus, assume henceforth that $f \geq 0$. For any $n \geq 0$, we have

$$(P^{2n}f)(X_0) := \mathbf{E}\left(f(X_{2n}) \mid X_0\right) = \mathbf{E}\left(\mathbf{E}\left(f(X_{2n}) \mid X_n\right) \mid X_0\right) = \mathbf{E}\left(R_n \mid X_0\right), \tag{9}$$

where $R_n := \mathbf{E}\left(f(X_{2n}) \mid X_n\right)$. Since $X_0 \sim \pi$, by reversibility,

$$(X_n, X_{n+1}, \ldots, X_{2n}) \quad \text{and} \quad (X_n, X_{n-1}, \ldots, X_0)$$

have the same law. Hence,

$$R_n = \mathbf{E}\left(f(X_{2n}) \mid X_n\right) = \mathbf{E}\left(f(X_0) \mid X_n\right) = \mathbf{E}\left(f(X_0) \mid X_n, X_{n+1}, \ldots\right), \tag{10}$$

where the second equality in (10) follows by the Markov property. Fix for now an integer $N \geq 0$. By (10), $(R_n)_{n=0}^N$ is a reverse martingale, i.e. $(R_{N-n})_{n=0}^N$ is a martingale. By Doob's $L^p$ maximal inequality (e.g. [1], Theorem 5.4.3),

$$\left\|\max_{0 \leq n \leq N} R_n\right\|_p \leq \frac{p}{p-1}\|R_0\|_p = \frac{p}{p-1}\|f(X_0)\|_p. \tag{11}$$

Define $h_N := \max_{0 \leq n \leq N} P^{2n}f$. By (9),

$$h_N(X_0) = \max_{0 \leq n \leq N} \mathbf{E}\left(R_n \mid X_0\right) \leq \mathbf{E}\left(\max_{0 \leq n \leq N} R_n \,\Big|\, X_0\right). \tag{12}$$

Recall that conditional expectation is an $L^p$ contraction (e.g. [1], Theorem 5.1.4). Thus, taking $L^p$ norms in (12) and applying (11), we have

$$\|h_N(X_0)\|_p \leq \left\|\max_{0 \leq n \leq N} R_n\right\|_p \leq \frac{p}{p-1}\|f(X_0)\|_p,$$

where in the first inequality we have implicitly used $f \geq 0$. Taking $N \to \infty$ and applying the monotone convergence theorem, we conclude that

$$\mathbf{E}_\pi\left|\sup_{n \geq 0} P^{2n}f\right|^p = \left\|\lim_{N \to \infty} h_N(X_0)\right\|_p^p \leq \left(\frac{p}{p-1}\right)^p \|f(X_0)\|_p^p = \left(\frac{p}{p-1}\right)^p \mathbf{E}_\pi |f|^p.$$

$\square$

## 9.2 Proof of Lemma 4.1

*Proof.* Let $Q$ be the transition matrix of $X$, and let $\pi$ be its stationary distribution. Define the inner product $\langle \cdot, \cdot \rangle_\pi$ on $\mathbb{R}^S$ by $\langle f, g \rangle_\pi = \sum_{s \in S} \pi(s) f(s) g(s)$. Note that by reversibility, we have

$$\langle Qf, g \rangle_\pi = \sum_{s \in S} \pi(s) \left(\sum_{s' \in S} Q(s, s') f(s')\right) g(s)$$

$$= \sum_{s \in S} \sum_{s' \in S} \pi(s') Q(s', s) f(s') g(s)$$

$$= \sum_{s' \in S} \pi(s') \left( \sum_{s \in S} Q(s', s) g(s) \right) f(s') = \langle f, Qg \rangle_\pi.$$

Hence, $Q$ is self-adjoint with respect to the $\pi$ inner product. Note that because $Q$ is stochastic, its eigenvalues lie in the interval $[-1, 1]$.

Now, let $\tilde{S} = S \setminus U$, and let $\tilde{Q}$ be the transition matrix of $X$ killed upon hitting $U$. That is, for $f \in \mathbb{R}^{\tilde{S}}$,

$$(\tilde{Q} f)(s) = \sum_{s' \in \tilde{S}} p(s, s') f(s').$$

If we regard $Q$ as a symmetric bilinear form, then $\tilde{Q}$ is the restriction of $Q$ onto the subspace $\mathbb{R}^{\tilde{S}}$. Hence, $\tilde{Q}$ is also a symmetric bilinear form, and if $Q$ is positive semidefinite, then so is $\tilde{Q}$.

Let $\tilde{X}_t$ be the walk started at $x$ and killed upon hitting $U$, and let $f_t(s) = \frac{1}{\pi(s)} \mathbf{P}(\tilde{X}_t = s)$ (we are guaranteed that $\pi(s) > 0$ by irreducibility). Note that $f_0(s) = \frac{1}{\pi(x)} \delta_{x,s}$, and

$$f_t(s) = \frac{1}{\pi(s)} \mathbf{P}(\tilde{X}_t = s) = \frac{1}{\pi(s)} \sum_{s' \in \tilde{S}} p(s', s) \mathbf{P}(\tilde{X}_{t-1} = s')$$

$$= \sum_{s' \in \tilde{S}} p(s, s') f_{t-1}(s') = (\tilde{Q} f_{t-1})(s),$$

so $f_t = \tilde{Q}^t f_0$. We thus have

$$\mathbf{P}(X_t = x, \tau(U) > t) = \mathbf{P}(\tilde{X}_t = x) = \pi(x) f_t(x)$$
$$= \pi(x) \cdot \langle f_t, f_0 \rangle_\pi = \pi(x) \cdot \langle \tilde{Q}^t f_0, f_0 \rangle_\pi.$$

Since $\tilde{Q}$ is self-adjoint with respect to $\langle \cdot, \cdot \rangle_\pi$, it has an orthonormal basis (in the $\pi$ inner product) of eigenvectors $\ell_1, \ldots, \ell_{|\tilde{S}|} \in \mathbb{R}^{\tilde{S}}$. For each $i$, let $\lambda_i$ be the eigenvalue corresponding to $\ell_i$. We can write $f_0 = a_1 \ell_1 + \cdots + a_{|\tilde{S}|} \ell_{|\tilde{S}|}$, and the last expression becomes

$$\pi(x) \sum_{i=1}^{|\tilde{S}|} a_i^2 \lambda_i^t.$$

This proves the lemma. $\square$